\documentclass[12pt]{amsart}
\usepackage{amsmath}
\usepackage{amsxtra}
\usepackage{amscd}
\usepackage{amsthm}
\usepackage{amsfonts}
\usepackage{amssymb}
\usepackage{eucal}
\usepackage{mathrsfs}
\usepackage{hyperref}
\usepackage{pstricks,pst-plot}

\def\C{{\mathbb C}}

\def\P{{\mathbb P}}

\def\Z{{\mathbb Z}}

\newtheorem{theorem}{Theorem}[section]
\newtheorem{definition}{Definition}[section]
\newtheorem{lemma}[theorem]{Lemma}
\newtheorem{remark}{Remark}

\newtheorem{corollary}[theorem]{Corollary}

\title{On higher syzygies of ruled surfaces}
\author{Euisung Park }
\address {Euisung Park : School of Mathematics, Korea Institute for Advanced Study,
270-43 Cheongryangri-dong, Dongdaemun-gu, Seoul 130-012, Republic
of Korea,} \email{puserdos@kias.re.kr}
\thanks{The author was supported by Korea Research Foundation Grant (KRF-2002-070-C00003).}

\begin{document}

\thispagestyle{empty} \maketitle

\setcounter{page}{1}

\begin{abstract}
We study higher syzygies of a ruled surface $X$ over a curve of
genus $g$ with the numerical invariant $e$. Let $L \in
\mbox{Pic}X$ be a line bundle in the numerical class of $aC_0
+bf$. We prove that for $0 \leq e \leq g-3$, $L$ satisfies
property $N_p$ if $a \geq p+2$ and $b-ae \geq 3g-1-e+p$ and for $e
\geq g-2$, $L$ satisfies property $N_p$ if $a \geq p+2$ and
$b-ae\geq 2g+1+p$. By using these facts, we obtain Mukai type
results. For ample line bundles $A_i$, we show that $K_X + A_1 +
\cdots + A_q$ satisfies property $N_p$ when $0 \leq e <
\frac{g-3}{2}$ and $q \geq g-2e+1 +p$ or when $e \geq
\frac{g-3}{2}$ and $q \geq p+4$. Therefore we prove Mukai's
conjecture for ruled surface with $e \geq \frac{g-3}{2}$. Also we
prove that when $X$ is an elliptic ruled surface with $e \geq 0$,
$L$ satisfies property $N_p$ if and only if $a \geq 1$ and
$b-ae\geq 3+p$.
\end{abstract}

\tableofcontents

\section{Introduction}
\noindent \noindent In this article we study higher syzygies of
irrational ruled surfaces. Let $X$ be a smooth projective variety
and $L \in \mbox{Pic}X$ a very ample line bundle. Consider the
embedding
\begin{equation*}
X \hookrightarrow \P H^0 (X,L) = \P
\end{equation*}
defined by the complete linear system of $L$. Let $S$ be the
homogeneous coordinate ring of $\P$ and consider the graded
$S$-module $E=\oplus_{n \in \Z} H^0 (X,L^n)$. Let
\begin{eqnarray*}
\cdots \rightarrow \oplus_j S^{k_{i,j}}(-i-j) \rightarrow  \cdots
\rightarrow \oplus_j S^{k_{1,j}}(-1-j) \rightarrow \oplus_j
S^{k_{0,j}}(-j) \rightarrow E \rightarrow 0
\end{eqnarray*}
be a minimal graded free resolution of $E$. A main goal to study
higher syzygies is to interpret the information carried by the
graded Betti number $k_{i,j}'s$. In particular, the distribution
of zeroes in the Betti table enables us to understand the geometry
of $X \hookrightarrow \P$. Along this point of view, many people
have studied the so-called property $N_p$ which means that first
few modules of syzygies are as simple as possible.\\

\begin{definition}[Green-Lazarsfeld, \cite{GL1}]
$(X,L)$ satisfies property $N_p$ if  $k_{i,j} = 0$ for $0 \leq i
\leq p$ and  $j \geq 2$. Equivalently, property $N_p$ holds for $X
\hookrightarrow \P$ if $E$ admits a minimal free resolution of the
form
\begin{eqnarray*}
\cdots \rightarrow S^{m_p}(-p-1) \rightarrow \cdots \rightarrow
S^{m_2}(-3) \rightarrow  S^{m_1}(-2) \rightarrow S \rightarrow E
\rightarrow 0.\\
\end{eqnarray*}
\end{definition}

\noindent Therefore property $N_0$ holds if and only if $X
\hookrightarrow \P H^0 (X,L)$ is a projectively normal embedding,
property $N_1$ holds if and only if property $N_0$ is satisfied
and the homogeneous ideal is generated by quadrics, and property
$N_p$ holds for $p \geq 2$ if and only if it has property $N_0$
and $N_1$ and the $k^{th}$ syzygies among the quadrics are
generated by linear syzygies for all $1 \leq k \leq p-1$.\\

Concerning higher syzygies, M. Green\cite{Green} obtained the
first general result. He proved that if $C$ is a smooth curve of
genus $g$ and if $\mbox{deg}(L) \geq 2g+1+p$, then $(C,L)$
satisfies property $N_p$. This result was rediscovered by M. Green
and R. Lazarsfeld\cite{GL2}. Also they classified all pairs
$(C,L)$ for which $``2g+1+p"$  theorem is optimal.

\begin{theorem}[Green-Lazarsfeld,
\cite{GL2}]\label{thm:classification} Let $L$ be a line bundle of
degree $2g+p~(p\geq1)$ on a smooth projective curve $X$ of genus
$g$, defining an embedding $C \hookrightarrow \P H^0
(C,L)=\P^{g+p}$. Then $(C,L)$ fails to satisfy property $N_p$ if
and only if either
\begin{enumerate}
\item[$(i)$] $C$ is hyperelliptic or \item[$(ii)$] $C
\hookrightarrow \P^{g+p}$ has a $(p+2)$-secant $p$-plane, i.e.,
$H^0 (C,L-K_C) \neq 0$ or equivalently $L=K_C + D$ for some
effective divisor $D$ of degree $p+2$.
\end{enumerate}
\end{theorem}

\noindent Then Mukai observed that Green's Theorem implies for
ample line bundle $A$ on a smooth curve $C$, $K_C + (p+3)A$
satisfies property $N_p$.
And he has conjectured that\\
\begin{enumerate}
\item[$(\star)$] {\bf Mukai's Conjecture.} For a smooth projective
surface $S$ and an ample line bundle $A
\in \mbox{Pic}S$,  $K_S + (4+p)A$ satisfies property $N_p$.\\
\end{enumerate}
Though this conjecture is still open, some progress has been made
by D. Butler for ruled surfaces\cite{Butler} and  by F. J. Gallego
and B. P. Purnaprajna for elliptic ruled
surfaces\cite{GP1}\cite{GP2}, surfaces of nonnegative Kodaira
dimension\cite{GP3}\cite{GP5} and rational surfaces\cite{GP4}.\\

The aim of this article is to study higher syzygies of ruled
surfaces over an irrational curve with numerical invariant $e \geq
0$. More precisely, we refine results in \cite{Butler} and
\cite{GP2}. We will follow the notation and terminology of R.
Hartshorne's book \cite{H}, V $\S 2$. Let $C$ be a smooth
projective curve of genus $g$ and let $\mathcal{E}$ be a vector
bundle of rank $2$ on $C$ which is normalized, i.e., $H^0
(C,\mathcal{E}) \neq 0$ while $H^0 (C,\mathcal{E} \otimes
\mathcal{O}_C (D))=0$ for every divisor $D$ of negative degree. We
set
\begin{equation*}
\mathfrak{e}=\wedge^2 \mathcal{E}~~~~\mbox{and}~~~~e = -
\mbox{deg}(\mathfrak{e}).
\end{equation*}
Let $X = \P_C (\mathcal{E})$ be the associated ruled surface with
projection morphism $\pi : X \rightarrow C$. We fix a minimal
section $C_0$ such that $\mathcal{O}_X (C_0)=\mathcal{O}_{\P_C
(\mathcal{E})} (1)$. For $\mathfrak{b} \in \mbox{Pic}C$,
$\mathfrak{b}f$ denote the pullback of $\mathfrak{b}$ by $\pi$.
Thus any element of $\mbox{Pic}X$ can be written
$aC_0+\mathfrak{b}f$ with $a\in \Z$ and $\mathfrak{b} \in
\mbox{Pic}C$ and any element of $\mbox{Num}X$ can be written $aC_0
+bf$ with $a,b \in \Z$.

When $g=1$, it is proved by Yuko Homma\cite{Homma1}\cite{Homma2}
and Gallego-Purnaprajna\cite{GP1} that property $N_0$ and $N_1$
are characterized in terms of the intersection number of $L$ with
a minimal section, a fiber and the anticanonical curve. See Remark
1.2 and 1.3.

\begin{theorem}\label{thm:genpres}
Let $X$ be an elliptic ruled surface. Let $L \in \mbox{Pic}X$ be a
line bundle in the numerical class of $aC_0 +bf$.\\
$(1)$ {\bf(Yuko Homma, \cite{Homma1}\cite{Homma2})} If $e \geq 0$,
then $L$ is normally generated if and only if

$a\geq1$ and $b-ae \geq 3$. If $e=-1$, then $L$ is normally
generated if and only if $a\geq1$,

$a+b \geq 3$ and $a+2b \geq 3$.\\
$(2)$ {\bf (F. J. Gallego and  B. P. Purnaprajna, \cite{GP1})} If
$e \geq 0$, then $L$ satisfies property

$N_1$ if and only if $a\geq1$ and $b-ae \geq 4$. If $e=-1$, then
$L$ satisfies property $N_1$ if

and only if $a\geq1$, $a+b \geq 4$ and $a+2b \geq 4$.
\end{theorem}

\noindent For higher syzygies, F. J. Gallego and B. P.
Purnaprajna\cite{GP2} obtained the following:
\begin{theorem}[F. J. Gallego and B. P. Purnaprajna, \cite{GP2}]\label{thm:GP}
Let $X$ be an elliptic ruled surface with the numerical invariant $e$.\\
$(1)$ Let $L$ be a line bundle in the numerical class of $aC_0
+bf$.
\begin{enumerate}
\item[$(a)$] If $e=-1$ and $a \geq p+1,~a+b \geq 2p+2$, and $a+2b
\geq 2p+2$, then $L$ satisfies property $N_p$. \item[$(b)$] If
$e\geq0$ and $a \geq p+1,~b-ae \geq 2p+2$, then $L$ satisfies
property $N_p$.
\end{enumerate}
$(2)$ If $B_1, \cdots,B_{p+1} \in \mbox{Pic}X$ are ample and base
point free line bundles, then property

$N_p$ holds for $B_1 + \cdots + B_{p+1}$.\\
$(3)$ If $A_1, \cdots,A_{2p+3} \in \mbox{Pic}X$ are ample line
bundles, then $K_X + A_1 + \cdots + A_{2p+3}$

satisfies property $N_p$.
\end{theorem}

\noindent Also they conjectured the following:\\

\noindent {\bf Conjecture.} (F. J. Gallego and B. P. Purnaprajna,
\cite{GP2}) Let $X$ be an elliptic ruled surface and $L \in
\mbox{Pic}X$ a line bundle in the numerical class $aC_0 + bf$.

$(1)$ If $e \geq 0$, then $L$ satisfies property $N_p$ if and only
if $a \geq 1$ and $b-ae \geq 3+p$.

$(2)$ If $e=-1$, then $L$ satisfies property $N_p$ if and only if
$a \geq 1$, $a+b \geq 3+p$, and

$a+2b \geq p+3$. \\

\noindent {\bf Remark 1.1.} When $e \geq 0$, $\mbox{deg}
(L|_{C_0}) = L \cdot C_0 = b-ae$. And $(C,L|_{C_0})$ satisfies
property $N_p$ if and only if $\mbox{deg} (L|_{C_0}) \geq 3+p$ by
Theorem \ref{thm:classification}. Therefore this conjecture
suggests that $(X,L)$ satisfies property $N_p$ if and
only if $(C,L|_{C_0})$ satisfies property $N_p$.\\

\noindent {\bf Remark 1.2.} When $e = -1$, note that there exists
a smooth elliptic curve $E \subset X$ such that $E \equiv 2C_0
-f$. See Proposition 3.2 in \cite{GP1}. Also $\mbox{deg}
(L|_{C_0}) = L \cdot C_0 = a+b$ and $\mbox{deg} (L|_{E}) = L \cdot
E = a+2b$. Thus $(C,L|_{C_0})$ satisfies property $N_p$ if and
only if $a+b \geq 3+p$, and $(E,L|_E)$ satisfies property $N_p$ if
and only if $a+2b \geq 3+p$ by Theorem \ref{thm:classification}.
Therefore this conjecture suggests that when $e =-1$, $(X,L)$
satisfies property $N_p$ if and only if $(C,L|_{C_0})$ and $(E,L|_{E})$ satisfy property $N_p$. \\

\noindent This conjecture has been solved for $p=0,1$ by Theorem
\ref{thm:genpres}. And our first main result is that this
conjecture is true when $e \geq 0$ or $e=-1$ and $a=1$:

\begin{theorem}\label{thm:elliptic}
Let $X$ be an elliptic ruled surface and $L \in \mbox{Pic}X$ a
line bundle in the numerical class $aC_0 + bf$.\\
$(1)$ If $e\geq 0$, then $L$ satisfies property $N_p$ if and only if $a\geq1$ and $b-ae\geq3+p$.\\
$(2)$ If $e=-1$, then $L$ satisfies property $N_p$ if $a\geq1$ and
$a+2b\geq 5+2p$.\\
$(3)$ If $e=-1$ and $a=1$, then $L$ satisfies property $N_p$ if
and only if $b\geq 2+p$.
\end{theorem}

We remove the condition $a \geq p+1$ in Theorem \ref{thm:GP}.(1).
Theorem \ref{thm:elliptic} says that if $e \geq 0$ or if $e=-1$
and $a=1$, then $(X,L)$ satisfies property $N_p$ if and only if
$(C,L_C)$ satisfies property $N_p$ where $L_C$ is the restriction
of $L$ to $C_0$. In particular, property $N_p$ is characterized in
terms of the intersection number of $L$ with a minimal section and
a fiber. Also we classify all projectively normal elliptic surface
scrolls which satisfies property $N_p$.

\begin{corollary}\label{cor:amplecaseelliptic}
Let $L=K_X + A_1 + \cdots +A_q$ where $A_i$ is ample and $q \geq
3$. If $e \geq 0$, then $L$ satisfies property $N_{p}$ if $q \geq
3-e+p$.
\end{corollary}

\begin{corollary}\label{cor:bpfcase}
Let $L=B_1 \otimes \cdots B_{p+1}$ be a line bundle on $X$ such
that each $B_i$ is ample and base point free and $p \geq 1$. If $e
\geq 0$, then $L$ satisfies property $N_{2p-1}$.
\end{corollary}

Corollary \ref{cor:amplecaseelliptic} and \ref{cor:bpfcase}
refines numerical bounds in Theorem \ref{thm:GP}.(2) and (3) when
$e \geq 0$. They are obtained by using Theorem \ref{thm:elliptic}
since all base point free line bundles and ample line bundles on
$X$ are classified. Corollary \ref{cor:amplecaseelliptic}.$(1)$
shows that Mukai's conjecture is true of elliptic ruled surfaces
with $e \geq 0$. Also this is optimal. Indeed let $A=aC_0 + bf$ be
an ample line bundle such that $b-ae=1$. Then $K_X + (3+p-e)A$
fails to
satisfy property $N_{p+1}$ by Theorem \ref{thm:elliptic}.\\

\noindent {\bf Remark 1.3.} Assume that $e=-1$ and let $A=aC_0 +
bf$ be an ample line bundle. If $a+2b \geq 2$, i.e., $A$ is base
point free, then $K_X + qA$ satisfies property $N_p$ for $q \geq
3+p$ since $K_X + qA = (aq-2)C_0 + (bq+1)f$. Therefore for
elliptic ruled surface with $e=-1$, Mukai's conjecture is true for
base point free ample line bundles, and the only remaining case is
when $a+2b=1$, that is, $A$ is ample but not base point
free$($e.g. $D=C_0)$.\\

\noindent \noindent When $C$ is a curve of genus $g \geq 2$, the
author has obtained that

\begin{theorem}[Corollary 4.4, \cite{ES}]\label{thm:EShighergenus}
Let $X$ be a ruled surface over a curve $C$ of genus $g \geq 2$
and $L \in \mbox{Pic}X$ a line bundle in the numerical class of
$aC_0 +bf$.\\
$(1)$ When $e \geq 0$, $L$ satisfies property $N_p$ if $b-ae \geq
3g-1+p$.\\
$(2)$ When $e < 0$, $L$ satisfies property $N_p$ if $2b-ae \geq
6g-2+2p$.\\
\end{theorem}

\noindent And in this article we prove the following sharper
result for $e \geq 0$:

\begin{theorem}\label{thm:main}
Let $X$ be a ruled surface over a curve $C$ of genus $g \geq 2$.
Let $L \in \mbox{Pic}X$ be a line bundle in the numerical class of
$aC_0 +bf$ and $L_C$ the restriction of $L$ to $C_0$.\\
$(1)$ When $0 \leq e \leq g-3$, $L$ satisfies property $N_p$ if $a \geq p+2$ and $b-ae\geq 3g-1-e+p$.\\
$(2)$ When $e \geq g-2$, $L$ satisfies property $N_p$ if $a \geq p+2$ and $b-ae\geq 2g+1+p$.\\
$(3)$ When $e \geq g-1$, assume that $a \geq p+3$ and $b-ae =
2g+p$. Then $(X,L)$

satisfies property $N_p$ if and only if $(C,L_C)$ satisfies
property $N_p$.
\end{theorem}

\noindent Note that $\mbox{deg}(L_C) = b-ae$. Therefore Theorem
\ref{thm:classification} and Theorem
\ref{thm:main}.$(3)$ guarantees the following:\\\\
$(*)$ Let $L \in \mbox{Pic}X$ be a line bundle in the numerical
class of $aC_0 +bf$ and let $L_C$

be the restriction of $L$ to $C_0$. When $e \geq g-1$ and $a \geq
p+2$, if $b-ae = 2g+p$

then $(X,L)$ fails to satisfy property $N_p$ if and only if either
\begin{enumerate}
\item[$(i)$] $C$ is hyperelliptic or
\item[$(ii)$] $L_C =K_C + D$ for some effective divisor $D$ of degree $p+2$.\\
\end{enumerate}

\noindent {\bf Example 1.1.} Assume that $e \geq g-1$ and fix a
positive integer $p$. Let $L=aC_0 + \mathfrak{b}f \in \mbox{Pic}X$
be such that $a \geq p+3$ and
\begin{equation*}
\mathfrak{b} =K_C - a \mathfrak{e} + D
\end{equation*}
and $D$ is an effective divisor of degree $2+p$. Then by Theorem
\ref{thm:classification}, $L_C = K_C +D$ does not satisfy Property
$N_p$. So Theorem \ref{thm:main} implies that $(X,L)$ satisfies
property $N_{p-1}$ and property $N_p$ does not hold.\\

\noindent Using Theorem \ref{thm:main} we prove the following
Mukai type result.

\begin{corollary}\label{cor:amplecase}
Let $X$ be a ruled surface over a curve $C$ of genus $g \geq 2$
with $e \geq 0$. Let $L=K_X + A_1 + \cdots +A_q$ where $A_i$ is
ample.

$(1)$ When $0 \leq e \leq g-3$, $L$ satisfies property $N_p$ if $q
\geq \mbox{max} \{4,g+1-2e \} +p$.

$(2)$ When $e \geq g-2$, $L$ satisfies property $N_p$ if $q \geq
4+p$. \\
In particular, Mukai's conjecture holds for $X$ if $e \geq
\frac{g-3}{2}$.
\end{corollary}

\noindent {\bf Remark 1.4.} Concerning Mukai's conjecture for
ruled surfaces with arbitrary $e$, D. Butler's result\cite{Butler}
says that $K_X + (4+4p)A$ satisfies property $N_p$. Also for ruled
surfaces with $e \geq 0$, his work says that $aC_0 + bf$ satisfies
property $N_p$ if $a \geq p+1$ and $b-ae \geq 2g+2p$.\\

\noindent {\bf Remark 1.5.} Let $X$ be a ruled surface over a
curve of genus $g$. Corollary \ref{cor:amplecase} shows that if
$g=2,3$ and $e \geq 0$, then Mukai's conjecture holds for $X$. \\

To study the distribution of zeros in the Betti table, we use the
so-called Koszul cohomology developed by Mark Green\cite{Green}.
In particular, this method enables us to show some vanishing of
Betti numbers by the vanishing of cohomology groups of certain
vector bundles. We reduce the problem to show vanishing of
cohomology groups on ruled surfaces to that on curves. Then we use
some numerical conditions to kill higher cohomology groups of
vector bundles on curves.

It seems the most interesting part of our result that when $e \geq
0$, higher syzygies of ruled surfaces is closely related to that
of the minimal section $C_0$. This follow from the long exact
sequence of cohomology groups induced by the short exact sequence
\begin{equation*}
0 \rightarrow \mathcal{O}_X (-C_0) \rightarrow \mathcal{O}_X
\rightarrow \mathcal{O}_C \rightarrow 0,
\end{equation*}
associated to the minimal section $C_0$ on $X$. On a variety $Z$
and a line bundle $L \in \mbox{Pic}Z$, a curve $C \subset Z$ is
said to be \textit{extremal} with respect to $L$ and property
$N_p$ if $(X,L)$ and $(C,L_C)$ satisfy property $N_p$ but not
property $N_{p+1}$. See Remark 1.5 in \cite{GP4}. When $X$ is a
rational surface with anti-canonical divisor nef, F. J. Gallego
and B. P. Purnaprajna\cite{GP4} proved for any very ample $L \in
\mbox{Pic}X$, a smooth curve $C$ in $|L|$ is extremal with respect
to $L$ by showing that the embedding of $X$ under $L$ is
arithmetically Cohen Macaulay. Applying this to rational ruled
surfaces, they obtained the following:

\begin{theorem}[F. J. Gallego and B. P. Purnaprajna,
\cite{GP4}]\label{thm:GPrational} Let $X$ be the rational ruled
surface associated to $\mathcal{O}_{\P^1} \oplus
\mathcal{O}_{\P^1} (-e)$. Let $L$ be a line bundle in the
numerical class of $aC_0 +bf$.
\begin{enumerate}
\item[$(1)$] If $a=1$ or $e=0$ and $b=1$, then $L$ satisfies
property $N_p$ for all $p \geq 0$.
\item[$(2)$] If $e=0$ and $a,b \geq 2$ or $e \geq 1$ and $a
\geq2$, then $L$ satisfies property $N_p$ if and only if $2a+2b-ae
\geq 3+p$.
\end{enumerate}
\end{theorem}

\noindent Unfortunately this does not hold for ruled surfaces over
an irrational curve. More precisely, $h^1 (X,\mathcal{O}_X)=g \geq
1$ and hence $(X,L)$ cannot be arithmetically Cohen Macaulay. So
we prove a similar statement. And our results implies that for $e
\geq 0$ the minimal section $C_0$ is the extremal curve. Obviously
the existence of extremal curves enables us to apply fruitful
results about higher syzygies of curves. When $e < 0$, there may
be more extremal curves. For an example, on an elliptic ruled
surface $X$ with $e=-1$ there exists a smooth elliptic curve $E
\subset X$ such that $E \equiv 2C_0 -f$. And Theorem
\ref{thm:genpres} and Theorem \ref{thm:elliptic}.$(3)$ make it
affirmative that $C_0$ and $E$ are the extremal curves with
respect to $L$ as conjectured by F.
J. Gallego and B. P. Purnaprajna.\\

The organization of this paper is as follows. In $\S 2$, we review
some necessary elementary facts to study higher syzygies of ruled
surfaces. $\S 3$ is devoted to give numerical criteria for
Property $N_p$. In $\S 4$ we suggest some open questions related
to our results.\\

\section{Preliminary}
\subsection{Notations and Conventions}
Throughout this paper the following is assumed.\\

\noindent $(1)$ All varieties are defined over the complex number
field $\C$.\\

\noindent $(2)$ For a finite dimensional $\C$-vector space $V$,
$\P(V)$ is the projective space of one-

dimensional quotients of $V$.\\

\noindent $(3)$ When a variety $X$ is embedded in a projective
space, we always assume that

it is non-degenerate, i.e. it does not lie
in any hyperplane.\\

\noindent $(4)$ When a projective variety $X$ is embedded in a
projective space $\P^r$ by a very ample

line bundle $L \in \mbox{Pic} X $, we may write $\mathcal{O}_X
(1)$ instead of $L$ as long as there is no

confusion.\\

\subsection{The slope of vector bundles on a curve}
Let $C$ be a smooth projective curve of genus $g$. For a vector
bundle $\mathcal{F}$ on $C$, the slope $\mu(\mathcal{F})$ is
defined by $\deg(\mathcal{F}) / \mbox{rank}(\mathcal{F})$. Also
the maximal slope and the minimal slope are defined as follows:
\begin{equation*}
\mu^+ (\mathcal{F})= \mbox{max}\{ \mu(S) | 0 \rightarrow S
\rightarrow \mathcal{F}\} ~~~\mbox{and}~~~ \mu^-
(\mathcal{F})=\mbox{min}\{ \mu(Q) | \mathcal{F}\rightarrow Q
\rightarrow 0\}
\end{equation*}
$\mathcal{F}$ is called semistable if $\mu(\mathcal{F})=\mu^-
(\mathcal{F})$ and stable if $\mu(\mathcal{F})<\mu^-
(\mathcal{F})$. These notions  satisfy the following properties.

\begin{lemma}\label{lem:folklore}
For vector bundles $\mathcal{E}$, $\mathcal{F}$ and $\mathcal{G}$ on $C$,\\
$(1)$ $\mu^+ ( \mathcal{E} \otimes \mathcal{F}) = \mu^+
(\mathcal{E} ) +\mu^+ ( \mathcal{F})$.\\
$(2)$ $\mu^- ( \mathcal{E} \otimes \mathcal{F}) = \mu^-
(\mathcal{E} ) +\mu^- ( \mathcal{F})$.\\
$(3)$ $\mu^+ (S^\ell (\mathcal{E})) = \ell \mu^+ (\mathcal{E})$.\\
$(4)$ $\mu^- (S^\ell (\mathcal{E})) = \ell \mu^- (\mathcal{E})$.\\
$(5)$ $\mu^- (\wedge^\ell \mathcal{E}) \geq \ell \mu^-
(\mathcal{E})$.\\
$(6)$ If $\mu^- (\mathcal{E}) > 2g-2$, then $h^1 (C, \mathcal{E})=0$.\\
$(7)$ If $\mu^- (\mathcal{E})  > 2g-1$, then $\mathcal{E}$ is globally generated.\\
$(8)$ If $\mu^- (\mathcal{E})  > 2g$, then
$\mathcal{O}_{\P(\mathcal{E})} (1)$ is very ample.\\
$(9)$ If $0 \rightarrow \mathcal{E} \rightarrow \mathcal{F}
\rightarrow \mathcal{G} \rightarrow 0$ is an exact sequence, then
\begin{eqnarray*}
\mu^- (\mathcal{G}) \geq \mu^- (\mathcal{F}) \geq \mbox{min}\{
\mu^- (\mathcal{E}), \mu^- (\mathcal{G}) \}.\\
\end{eqnarray*}
\end{lemma}

\begin{proof}
See $\S 1$ and $\S 2$ in \cite{Butler}.
\end{proof}

\noindent Therefore for a vector bundle $\mathcal{E}$, if $\mu^-
(\mathcal{E}) > 2g-1$, then the evaluation map determines an exact
sequence of bundles:
\begin{equation*}
0 \rightarrow \mathcal{M}_{\mathcal{E}} \rightarrow H^0
(C,\mathcal{E}) \otimes \mathcal{O}_C \rightarrow \mathcal{E}
\rightarrow 0.
\end{equation*}

\noindent And Butler obtained the following very useful result:

\begin{theorem}[Butler, \cite{Butler}]\label{thm:Butlerestimation}
For a vector bundle $\mathcal{E}$ over $C$, if $\mu^-
(\mathcal{E})\geq 2g$, then $\mathcal{M}_{\mathcal{E}}$ satisfies
\begin{equation*}
\mu^- (\mathcal{M}_{\mathcal{E}}) \geq -\frac{\mu^-
(\mathcal{E})}{\mu^- (\mathcal{E})-g}.
\end{equation*}\\
\end{theorem}

\subsection{Regularity of vector bundles over ruled varieties }
Let $E$ be a vector bundle of rank $n+1$ over a smooth projective
variety $Y$ and let $X=\P_Y (E)$ with the projection map $\pi : X
\rightarrow Y$ and tautological line bundle $H$.

\begin{definition}
For a vector bundle $\mathcal{F}$ over $X$, we say that
$\mathcal{F}$ is $f$  $\pi$-regular when
\begin{equation*}
R^i \pi_* (\mathcal{F}(f-i)) = 0
\end{equation*}
for every $i \geq 1$.
\end{definition}

\noindent Here $\mathcal{F}(f-i) =\mathcal{F} \otimes H^{f-i}$. By
definition, a line bundle of the form $aH+\pi^* B$ is $(-a)$
$\pi$-regular. We present some basic facts about the
$\pi$-regularity.

\begin{lemma}[Lemma 3.2, \cite{Butler}]\label{lem:directimage} Let $\mathcal{F}$ and
$\mathcal{G}$ be two vector bundles on $X$
with $f$ and $g$ $\pi$-regularity, respectively.\\
$(1)$ $\mathcal{F} \otimes \mathcal{G}$ is $(f+g)$
$\pi$-regular.\\
$(2)$ If $f\leq 1$, then
\begin{equation*}
H^i (X,\mathcal{F})\cong H^i (Y,\pi_* \mathcal{F})\quad \mbox{for
all}\quad i\geq0.
\end{equation*}
$(3)$ If $f\leq0$ and $\widetilde{\mathcal{F}}=\pi^* (\pi_*
\mathcal{F})$, there is an exact sequence of vector bundles on $X$
\begin{equation*}
0 \rightarrow  \mathcal{K}_{\widetilde{\mathcal{F}}} \rightarrow
\widetilde{\mathcal{F}}  \rightarrow \mathcal{F} \rightarrow 0
\end{equation*}

where $\mathcal{K}_{\widetilde{\mathcal{F}}}$ is $1$ $\pi$-regular.\\
$(4)$ If $f \leq 0$ and $g \leq 0$, then there is a surjective map
\begin{equation*}
\pi_* \mathcal{F} \otimes \pi_*  \mathcal{G} \rightarrow \pi_*
(\mathcal{F} \otimes \mathcal{G}) \rightarrow 0.
\end{equation*}

In particular if $Y$ is a curve, $\mu^- (\pi_* (\mathcal{F}
\otimes \mathcal{G})) \geq \mu^- (\pi_* \mathcal{F}) +  \mu^-
(\pi_*  \mathcal{G})$.\\
\end{lemma}

\subsection{Ruled Surfaces}  Let $X$ be
a ruled surface over a smooth projective curve $C$ of genus $g$. We use notations in $\S 1$ and survey some basic facts about $X$:\\
\begin{enumerate}
\item[$(1)$] The restriction of $aC_0+\mathfrak{b}f$ to $C_0$
is $\mathfrak{b}+a\mathfrak{e}$.\\
\item[$(2)$] The canonical line bundle of $X$ is given by $K_X =
-2C_0 + (K_C +\mathfrak{e})f$.\\
\item[$(3)$] The arithmetic genus $p_a (X) = -g$ and the geometric genus $p_g (X) =0$.\\
\item[$(4)$] When $e \geq 0$, $\mu^- (\mathcal{E})=-e$. When $e
\leq -1$, $\mu^- (\mathcal{E})=-\frac{e}{2}$.\\
\item[$(5)$] $\mathcal{E}$ is nonstable if and only if $e >0$,
semistable if and only if $e \leq 0$ and stable if and only if $e
< 0$.\\
\item[$(6)$] Let $L \in \mbox{Pic}X$ be a line bundle in the
numerical class of $aC_0 +bf$.
    \begin{enumerate}
       \item[$(i)$] If $e \geq 0$, then $L$ is ample if and only
                       if $a \geq 1$ and $b -ae \geq 1$.
       \item[$(ii)$] If $e < 0$, then $L$ is ample if and only
                       if $a \geq 1$ and $2b -ae \geq 1$.\\
    \end{enumerate}
\end{enumerate}

\subsection{Koszul Cohomology} Let $V$ be a finite-dimensional complex vector space and
let $S(V)$ be the symmetric algebra on $V$. For a nonzero coherent
sheaf $\mathcal{F}$ on $\P=\P V$, consider the associated graded
$S$-module
\begin{equation*}
F=\oplus_{\ell \in \Z} H^0 (\P,\mathcal{F}(\ell))
\end{equation*}
and the minimal free resolution
\begin{eqnarray*}
\cdots \rightarrow  \oplus_{q \in \Z} S(V)(-q) \otimes M_{i,q}
\stackrel{\varphi_i}{\rightarrow}  \oplus_{q \in \Z} S(V)(-q)
\otimes M_{i-1,q} \rightarrow \cdots \\
\rightarrow \oplus_{q \in \Z} S(V)(-q) \otimes M_{1,q}
\stackrel{\varphi_1}{\rightarrow} \rightarrow \oplus_{q \in \Z}
S(V)(-q) \otimes M_{0,q} \stackrel{\varphi_0}{\rightarrow} F
\rightarrow 0
\end{eqnarray*}
of $F$. Put $k_{i,j} = \mbox{dim}_k M_{i,i+j}$. These integers are
called $Betti~~~~numbers$ of $F$. There is the following general
connection between syzygies and some cohomology groups.

\begin{theorem}[Theorem 4.5, \cite{Eisenbud}]\label{thm:exactsequence}
There is an exact sequence
\begin{equation*}
\begin{CD}
\wedge^{i+1}V \otimes H^0 ( \P, \mathcal{F} (j-1))
\stackrel{\alpha_{i,j}}{\rightarrow} H^0 (\P,\wedge^i \mathcal{M}
\otimes \mathcal{F} (j) )  \rightarrow \\
 H^1 (\P,\wedge^{i+1} \mathcal{M} \otimes \mathcal{F} (j-1) )
 \rightarrow  \wedge^{i+1} V \otimes H^1 ( \P, \mathcal{F} (j-1))
\end{CD}
\end{equation*}
with $\mbox{Coker}(\alpha_{i,j}) \cong M_{i,i+j}$  where
$\mathcal{M}=\Omega _{\P^r} (1)$ and $V=H^0 (\P^r ,
\mathcal{O}_{\P^r} (1))$. Therefore $k_{i,j}=\mbox{dim}_k
Coker(\alpha_{i,j})$ and we have the exact sequence
\begin{equation*}
0 \rightarrow M_{i,i+j} \rightarrow  H^1 (\P ,\wedge^{i+1}
\mathcal{M} \otimes \mathcal{F} (j-1) ) \rightarrow \wedge^{i+1} V
\otimes H^1 ( \P,\mathcal{F} (j-1)).
\end{equation*}
\end{theorem}

\noindent In this paper, we concern the case when $\mathcal{F}$ is
a coherent sheaf on a subvariety $X \subset \P$. Under this
situation, Theorem \ref{thm:exactsequence} guarantees that there
is an exact sequence
\begin{equation*}
\begin{CD}
\wedge^{i+1}V \otimes H^0 (X, \mathcal{F} (j-1))
\stackrel{\alpha_{i,j}}{\rightarrow} H^0 (X,\wedge^i \mathcal{M}
\otimes \mathcal{F} (j) )  \rightarrow \\
 H^1 (X,\wedge^{i+1} \mathcal{M} \otimes \mathcal{F} (j-1) )
 \rightarrow  \wedge^{i+1} V \otimes H^1 ( X, \mathcal{F} (j-1))
\end{CD}
\end{equation*}
with $\mbox{Coker}(\alpha_{i,j}) \cong M_{i,i+j}$.\\

\subsection{Cohomological interpretation of property $N_p$} We
review some cohomological criteria for property $N_p$. Let $X$ be
a smooth projective variety of dimension $n \geq 1$ and let $L \in
\mbox{Pic}X$ be a very ample line bundle. Consider the natural
short exact sequence
\begin{equation*}
0 \rightarrow \mathcal{M}_L \rightarrow H^0 (X,L) \otimes
\mathcal{O}_X \rightarrow L \rightarrow 0.
\end{equation*}

\begin{lemma}\label{lem:criterion}
Suppose that the ideal sheaf $~\mathcal{I}_{X/\P}~$ of $X
\hookrightarrow \P = \P H^0 (X,L)$ is $3$-regular in the sense of
Castelnuovo-Mumford, i.e., that $H^i (\P,\mathcal{I}_{X/\P}
(3-i))=0$ for all $i \geq 1$. Then for $p \leq
\mbox{codim}(X,\P)$, property $N_p$ holds for $L$ if and only if
$H^1 (X,\wedge^{p+1} \mathcal{M}_L \otimes L)=0$.
\end{lemma}

\begin{proof}
See $\S 1$ in \cite{GL2}.
\end{proof}

\begin{remark}\label{rmk:criterion}
Let $X$ be a smooth projective surface with geometric genus $0$,
i.e., $H^2 (X,\mathcal{O}_X)=0$. Let $L \in \mbox{Pic}X$ be a
normally generated very ample line bundle such that $H^1 (X,L)=0$.
Then  it is easy to check that $\mathcal{I}_{X/\P}$ is $3$-regular.\\
\end{remark}

\subsection{Higher syzygies of degenerate varieties}
Let $\Lambda \cong \P W \subset \P V$ be a linear subspace such
that $\mbox{codim}(W,V)=c$. It is easily checked that
\begin{equation*}
\Omega_{\P V} (1)|_{\P W} \cong \Omega_{\P W} (1) \oplus
\mathcal{O}_{\P W} ^c.
\end{equation*}
Now let $X \subset \P V$ be a smooth projective variety which is
indeed contained in $\Lambda$. Let the corresponding very ample
line bundle on $X$ be $L \in \mbox{Pic}X$. Consider the natural
short exact sequences
\begin{eqnarray*}
& 0 \rightarrow \mathcal{M}_V  \rightarrow V \otimes \mathcal{O}_X \rightarrow L \rightarrow  0&~~~~\mbox{and} \\
& 0 \rightarrow \mathcal{M}_W \rightarrow W \otimes \mathcal{O}_X
\rightarrow L \rightarrow  0.
\end{eqnarray*}
The above observation shows that $\mathcal{M}_V \cong
\mathcal{M}_W \oplus \mathcal{O}_{X} ^c$.

\begin{lemma}\label{lem:degenrate}
Under the situation just stated, assume that $H^1 (X,L^{\ell}) =0$
for all $\ell \geq 1$ and $W = H^0 (X,L)$. Then $(X,L)$ satisfies
property $N_p$ if and only if
\begin{equation*}
H^1 (X,\wedge^i \mathcal{M}_V \otimes L^{\ell}) = 0~~~~\mbox{for
$1 \leq i \leq p+1$ and $\ell \geq 1$.}
\end{equation*}
\end{lemma}

\begin{proof}
We use induction on $c$.\\
When $c=0$, this is a well-known criterion for Property $N_p$.\\
Assume that $c=1$. Since $\mathcal{M}_W = \mathcal{M}_V \oplus
\mathcal{O}_{X}$, $\wedge^i  \mathcal{M}_V \cong \wedge^i
\mathcal{M}_W \oplus \wedge^{i-1} \mathcal{M}_W$ which completes
the proof.\\
When $c \geq 2$, fix a filtration $H^0 (X,L)=W \subset W_1 \subset
\cdots \subset W_c = V$ of subspaces each having codimension one
in the next. Then assertion comes by repeating the process in the
case $c=1$.
\end{proof}

\section{Main Theorems}
\noindent Let $\pi : X \rightarrow C$ be the projection morphism.
Let $L = aC_0 +\mathfrak{b}f \in \mbox{Pic}X$ be a line bundle in
the numerical class of $aC_0 +bf$ such that $a \geq 1$ and $b+a
\mu^- * (\mathcal{E}) \geq 2g+1$. Let $\mathcal{F}$ denote $\pi_*
L = S^a (\mathcal{E})\otimes \mathfrak{b}$. Then $\mu^- (\pi_* L)
= b+a \mu^- * (\mathcal{E})$ and hence by Lemma
\ref{lem:folklore}, $L$ is very ample and $H^1 (X,L^{\ell})=0$ if
$\ell \neq 0$. Now consider the exact
sequences\\
\begin{eqnarray*}
0 \rightarrow  \mathcal{M}_{\mathcal{F}} \rightarrow  H^0
(C,\mathcal{F}) \otimes \mathcal{O}_C
\rightarrow  \mathcal{F} \rightarrow 0 &~~~~\mbox{and} \\
0 \rightarrow  \mathcal{M}_L \rightarrow H^0 (X,L) \otimes
\mathcal{O}_X \rightarrow L \rightarrow 0.
\end{eqnarray*}\\

\begin{lemma}\label{lem:generalvanishing}
Let $N \in \mbox{Pic}X$ be a line bundle in the numerical class of
$sC_0 +tf$. \\
$(1)$ For $m \geq 1$, $H^1 (X,\wedge^m \mathcal{M}_L \otimes N)=0$
if
\begin{equation*}
s \geq \mbox{min}\{m,a\}~~~~\mbox{and}~~~~t+s \mu^- (\mathcal{E})
> \frac{m(b+a\mu^- (\mathcal{E}))}{b+a\mu^- (\mathcal{E})-g}+2g-2.
\end{equation*}
$(2)$ For $m \geq 1$, $H^2 (X,\wedge^m \mathcal{M}_L \otimes N)=0$
if $s+1 \geq m$.
\end{lemma}

\begin{proof}
Consider the exact sequence
\begin{eqnarray*}
0 \rightarrow  \mathcal{K}_L \rightarrow \pi^* \mathcal{F}
\rightarrow L \rightarrow 0
\end{eqnarray*}
where $\mathcal{K}_L$ is a vector bundle of rank $a$ on $X$ which
is $1$ $\pi$-regular by Lemma \ref{lem:directimage}. Using Snake
Lemma , we have the following commutative diagram:
\begin{equation*}
\begin{CD}
&  &&  && 0 &\\
&  &&  && \downarrow &\\
& 0 && && \mathcal{K}_L &\\
& \downarrow &&  && \downarrow &\\
0 \rightarrow &\pi^* \mathcal{M}_\mathcal{F}& \rightarrow & H^0 (C,\mathcal{F}) \otimes \mathcal{O}_X & \rightarrow  & \pi^* \mathcal{F} & \rightarrow 0 \\
& \downarrow && \parallel \wr && \downarrow &\\
0 \rightarrow & \mathcal{M}_L & \rightarrow & H^0 (X,L) \otimes \mathcal{O}_X & \rightarrow & L & \rightarrow 0 \\
& \downarrow &&  && \downarrow &\\
& \mathcal{K}_L && && 0  &\\
& \downarrow  &&  && &\\
& 0 &&  && &.\\
\end{CD}
\end{equation*}\\
$(1)$ From the exact sequence
\begin{equation*}
0 \rightarrow \pi^* \mathcal{M}_{\mathcal{F}} \rightarrow
\mathcal{M}_L \rightarrow \mathcal{K}_L \rightarrow 0,
\end{equation*}
the desired vanishing is obtained if
\begin{equation*}
H^1 (X,\wedge^{m-i} \pi^* \mathcal{M}_{\mathcal{F}} \otimes
\wedge^i \mathcal{K}_L \otimes N)=0~~~~\mbox{for all $0 \leq i
\leq \mbox{min}\{m,a\}$}.
\end{equation*}
 Since we are in characteristic zero,
$\wedge^i \mathcal{K}_L$ is a direct summand of $T^i
\mathcal{K}_L$, and therefore it suffices to show that
\begin{equation*}
H^1 (X,\wedge^{m-i} \pi^* \mathcal{M}_{\mathcal{F}} \otimes T^i
\mathcal{K}_L \otimes N)=0~~~~\mbox{for all $0 \leq i \leq
\mbox{min}\{m,a\}$}.
\end{equation*}
Lemma \ref{lem:directimage} guarantees that $T^i \mathcal{K}_L$ is
$i$ $\pi$-regular and hence
\begin{equation*}
H^1 (X,\wedge^{m-i} \pi^* \mathcal{M}_{\mathcal{F}} \otimes T^i
\mathcal{K}_L \otimes N)=H^1 (C, \wedge^{m-i}
\mathcal{M}_{\mathcal{F}} \otimes \pi_* T^i \mathcal{K}_L \otimes
N)
\end{equation*}
from the assumption $s \geq \mbox{min}\{m,a\}$. So by Lemma
\ref{lem:folklore}, it suffices to show that
\begin{equation*}
\mu^- (\wedge^{m-i}  \mathcal{M}_{\mathcal{F}} \otimes \pi_* T^i
\mathcal{K}_L \otimes N)>2g-2.
\end{equation*}
First we claim that\\
\begin{enumerate}
\item[$(i)$] $\mu^- (\pi_* \mathcal{K}_L \otimes C_0) \geq b+a\mu^- (\mathcal{E})$.
\item[$(ii)$] $\mu^- (\pi_* T^i \mathcal{K}_L \otimes N) \geq i
(b+a\mu^- (\mathcal{E}))+t+s\mu^- (\mathcal{E})$.\\
\end{enumerate}
For $(i)$, apply Lemma 4.3 in \cite{Butler} to the case $V = L$
and $W = C_0$ in Butler's notation. Since $\mathcal{K}_L\otimes
C_0$ is $0$ $\pi$-regular and
\begin{equation*}
T^i \mathcal{K}_L \otimes N = T^i (\mathcal{K}_L \otimes C_0)
\otimes ((s-i)C_0 + tf),
\end{equation*}
repeated application of Lemma \ref{lem:directimage}.(4) shows that
\begin{eqnarray*}
\mu^- (\pi_* T^i \mathcal{K}_L \otimes N) & \geq & i \mu^- (\pi_*
\mathcal{K}_L \otimes C_0) + \mu^- (\pi_* \{(s-i)C_0 +tf \}) \\
& \geq & i (b+(a+1)\mu^- (\mathcal{E}))+ (s-i)\mu^- (\mathcal{E})+t \\
& = & i (b+a\mu^- (\mathcal{E}))+t+s\mu^- (\mathcal{E})
\end{eqnarray*}
which completes the proof of $(ii)$. From this claim,
\begin{eqnarray*}
\mu^- (\wedge^{m-i}  \mathcal{M}_{\mathcal{F}} \otimes \pi_* T^i
\mathcal{K}_L \otimes N)
   & \geq & (m-i) \mu^- (\mathcal{M}_{\mathcal{F}}) + i
(b+a\mu^- (\mathcal{E}))+t+s\mu^- (\mathcal{E})\\
   & \geq & m \mu^- (\mathcal{M}_{\mathcal{F}}) +t+s\mu^- (\mathcal{E}) \\
   & \geq & -m\frac{\mu^-
(\mathcal{F})}{\mu^- (\mathcal{F})-g}
+t+s\mu^-(\mathcal{E})~~\mbox{(Theorem
\ref{thm:Butlerestimation})}.
\end{eqnarray*}
Since  $\mu^- (\mathcal{F})= b +a\mu^- (\mathcal{E})$,
\begin{equation*}
- \frac{m(b+a\mu^- (\mathcal{E}))}{b+a\mu^- (\mathcal{E})-g}
+t+s\mu^- (\mathcal{E})
> 2g-2
\end{equation*}
implies the desired vanishing.\\
$(2)$ In the same way, $H^2 (X,\wedge^m \mathcal{M}_L \otimes
N)=0$ if
\begin{equation*}
H^2 (X,\wedge^{m-i} \pi^* \mathcal{M}_{\mathcal{F}} \otimes T^i
\mathcal{K}_L \otimes N)=0~~~~\mbox{for all $0 \leq i \leq m$}.
\end{equation*}
Since $T^i \mathcal{K}_L$ is $i$ $\pi$-regular, if $s+1 \geq m$,
then $T^i \mathcal{K}_L \otimes N$ is $1$ $\pi$-regular and hence
\begin{equation*}
H^2 (X,\wedge^{m-i} \pi^* \mathcal{M}_{\mathcal{F}} \otimes T^i
\mathcal{K}_L \otimes N)=H^2 (C,\wedge^{m-i}
\mathcal{M}_{\mathcal{F}} \otimes \pi_* T^i \mathcal{K}_L \otimes
N)=0.
\end{equation*}
by Lemma \ref{lem:directimage}. \qed \\\\

\noindent {\bf Proof of Theorem \ref{thm:elliptic}(1) and (2).} By
Lemma \ref{lem:generalvanishing}, $H^1 (X,\wedge^m \mathcal{M}_L
\otimes L)=0$ if
\begin{equation*}
b+a\mu^- (\mathcal{E})
>  \frac{m(b+a\mu^- (\mathcal{E}))}{b+a\mu^- (\mathcal{E})-1}
\end{equation*}
or equivalently $b+a\mu^- (\mathcal{E}) > m+1$. Therefore if
$b+a\mu^- (\mathcal{E}) > 2+p$, then $L$ satisfies property $N_p$
by Lemma \ref{lem:criterion} and Remark \ref{rmk:criterion}. Since
$\mu^- (\mathcal{E}) = -e$ if $e \geq 0$ and
$\mu^- (\mathcal{E}) = \frac{1}{2}$ if $e =-1$, it is proved that\\
\begin{enumerate}
\item[$(1)$] If $e\geq 0$, then $L$ satisfies property $N_p$ if
$a\geq1$ and $b-ae\geq3+p$.
\item[$(2)$] If $e=-1$, then $L$ satisfies property $N_p$ if
$a\geq1$ and $a+2b\geq5+2p$.\\
\end{enumerate}

\noindent To complete the proof it remains to show that when $e \geq 0$,\\
\begin{enumerate}
\item[$(*)$] if $b-ae = 3+p$ for some $p \geq 0$, then $L$
fails to satisfy Property $N_{p+1}$.\\
\end{enumerate}
Consider the exact sequence
\begin{equation*}
0 \rightarrow \mathcal{O}_X (-C_0) \rightarrow \mathcal{O}_X
\rightarrow \mathcal{O}_C \rightarrow 0.
\end{equation*}
Define graded $S$-modules $R_1 =\oplus_{\ell \in \Z} H^0
(X,\mathcal{O}_X (-C_0) \otimes L^{\ell})$,  $R_2 =\oplus_{\ell
\in \Z} H^0 (X,L^{\ell})$, and $R_3 =\oplus_{\ell \in \Z} H^0
(C,L_C ^{\ell})$ where $L_C$ is the restriction of $L$ to $C$. It
is easy to check that $H^1 (X,\mathcal{O}_X (-C_0) \otimes
L^{\ell})=0$ for all $\ell >0$(see Proposition 3.1 in \cite{GP1}).
So we have the exact sequence
\begin{equation*}
0 \rightarrow R_1 \rightarrow R_2 \rightarrow R_3 \rightarrow 0
\end{equation*}
of graded $S$-modules with maps preserving the gradings and hence
there is a long exact sequence
\begin{eqnarray*}
\cdots \rightarrow M_{1,q}(R_1,V) \rightarrow
M_{1,q}(R_2,V) \rightarrow M_{1,q}(R_3,V)\\
\rightarrow M_{0,q}(R_1,V) \rightarrow M_{0,q}(R_2,V) \rightarrow
M_{0,q}(R_3,V)\rightarrow 0.
\end{eqnarray*}
See Corollary (1.d.4) in \cite{Green}. We need the following part:
\begin{eqnarray*}
\cdots \rightarrow M_{p+1,p+3}(R_2,V) \rightarrow
M_{p+1,p+3}(R_3,V) \rightarrow M_{p,p+3}(R_1,V) \rightarrow \cdots
\end{eqnarray*}
Since $H^1 (X,\mathcal{O}_X (-C_0) \otimes L^2)=0$,
\begin{equation*}
M_{p,p+3}(R_1,V) =H^1 (X,\wedge^{p+1} \mathcal{M}_L \otimes
\mathcal{O}_X (-C_0) \otimes L^2)
\end{equation*}
by Theorem \ref{thm:exactsequence}. By Lemma
\ref{lem:generalvanishing},  $H^1 (X,\wedge^{p+1} \mathcal{M}_L
\otimes \mathcal{O}_X (-C_0) \otimes L^2)=0$ and hence
$M_{p,p+3}(R_1,V)=0$.

Assume that $L$ satisfies property $N_{p+1}$, then
$M_{p+1,p+3}(R_2,V)=0$ and hence $M_{p+1,p+3}(R_3,V)=0$ which
implies that $(C,L_C)$ satisfies property $N_{p+1}$. This is
impossible since $\mbox{deg} (L_C) =3+p$. So
$M_{p+1,p+3}(R_2,V)\neq0$ and hence property $N_{p+1}$ does not
hold for $(X,L)$. \qed\\

Now we prove Corollary \ref{cor:amplecaseelliptic} and
\ref{cor:bpfcase}. Recall numerical criteria for base point
freeness for line bundles on $X$.

\begin{lemma}\label{lem:bpf}
Let $L \in \mbox{Pic}X$ be a line bundle in the numerical class of
$aC_0 +bf$.\\
$(1)$ If $e\geq 0$, then $L$ is ample and base point free if and only if $a\geq1$ and $b-ae\geq2$.\\
$(2)$ If $e=-1$, then $L$ is ample and base point free if and only
if $a\geq1$, $a+b\geq2$ and

$a+2b \geq2$.
\end{lemma}

\begin{proof}
See Remark 3.5.3 in \cite{GP1}. \qed \\\\

\noindent {\bf Proof of Corollary \ref{cor:amplecaseelliptic}.}
Put $A_i \equiv a_i C_0 + b_i f$. Thus $a_i \geq 1$ and $b_i - a_i
e \geq 1$ by Lemma \ref{lem:bpf}. Also
\begin{equation*}
K_X + A_1 + \cdots +A_q \equiv (\sum a_i -2 ) C_0 + (\sum b_i
-e)f.
\end{equation*}
Since $(\sum b_i -e)- e(\sum a_i -2 )= \sum (b_i - a_i e) + e \geq
q+e$ for $e \geq 0$,
Theorem \ref{thm:elliptic} guarantees $(1)$.  \qed \\\\

\noindent {\bf Proof of Corollary \ref{cor:bpfcase}.} Put $B_i =
a_i C_0 + b_i f$. Thus $a_i \geq 1$ and $b_i - a_i e \geq 2$ by
Lemma \ref{lem:bpf}. Also
\begin{equation*}
L = (\sum a_i) C_0 + (\sum b_i) f~~~~\mbox{and}~~~~\sum b_i - e
\sum a_i  = \sum (b_i - e a_i) \geq 2p+2.
\end{equation*}
Then the assertion comes from Theorem \ref{thm:elliptic}.
\end{proof}

\noindent For a given normally generated very ample line bundle,
it is one of the most natural questions that for which $p$ the
line bundle satisfies property $N_p$. If $X$ is an elliptic ruled
surface with $e \geq 0$, this question is answered perfectly in
Theorem \ref{thm:elliptic}.(1). It seems very interesting that\\

\begin{enumerate}
\item[$(*)$] for a line bundle $L \in \mbox{Pic}X$ in the
numerical class of $aC_0 +bf$ with $a \geq 1$, property $N_p$
holds for $L$ if and only if property $N_p$ holds for the
restriction of $L$ to the minimal section $C_0$. \\
\end{enumerate}

\noindent That is, for elliptic ruled surfaces with $e\geq0$ the
higher syzygies of line bundles are closely connected with that of
the minimal section. Now we turn to the case $e=-1$. Recall that
if $X$ is an elliptic ruled surface with$e = -1$, then there
exists a smooth elliptic curve $E \subset X$ such that $E \equiv
2C_0 -f$.  It is proved by Yuko Homma\cite{Homma1}\cite{Homma2}
and Gallego-Purnaprajna\cite{GP1} that property $N_0$ and $N_1$
are characterized in terms of the intersection number of $L$ with
a minimal section, a fiber and the anticanonical curve $E$.
Theorem \ref{thm:genpres} makes the following conjecture affirmative.\\

\noindent {\bf Conjecture.} (F. J. Gallego and B. P. Purnaprajna,
\cite{GP2}) Let $X$ be an elliptic ruled surface with $e=-1$ and
let $L \in \mbox{Pic}X$ be a line bundle in the numerical class
$aC_0 + bf$. Then $L$ satisfies property $N_p$ if and only if $a
\geq 1$, $a+b \geq 3+p$, and $a+2b \geq p+3$. \\

\noindent For $p=2$, this is checked for $L$ in the numerical
class $C_0 +4f$ or $2C_0 + 3f$(\cite{GP2}, $\S 7$). Now we show
that this conjecture is true for scrolls, i.e., $a=1$.\\

\begin{figure}[t]
\begin{center}
\begin{pspicture}(-2,-2)(11,11)
\footnotesize \psset{xunit=.5cm,yunit=.5cm}

\psaxes[labels=none,ticks=none]{->}(0,0)(-3,-3)(20,17)

\uput[180](-.5,8){$p+\frac{5}{2}$} \uput[180](-.5,9){$p+3$}
\uput[180](-.5,14){$2p+2$} \uput[180](-.5,17){$f$}

\rput[t](2,-.5){1} \rput[t](6,-.5){$p+1$} \rput[t](9,-.5){$p+3$}
\rput[t](12,-.5){$2p-1~~$} \rput[t](14,-.5){$2p+2~~$}
\rput[t](16,-.5){$2p+5$} \rput(21,0){$C_0$}

\psline[linestyle=dashed](2,0)(2,7) \psline(2,7)(2,16)
\psline[linestyle=dashed](6,0)(6,8) \psline(6,8)(6,16)

\pspolygon[linestyle=none,fillstyle=vlines,hatchsep=6pt](6,16)(6,8)(14,0)(18,-2)(19,-2)(19,16)
\pspolygon[linestyle=none,fillstyle=hlines,hatchsep=6pt](2,16)(2,7)(19,-1.5)(19,16)
\psline[linestyle=dashed](12,0)(12,2)
\psline[arrows=*-*](0,9)(9,0) \psline[arrows=->](9,0)(13,-2)
\psline[arrows=*-*](0,14)(14,0) \psline[arrows=->](14,0)(18,-2)
\psline[arrows=*->](0,8)(19,-1.5)

\end{pspicture}
\end{center}
\caption{Elliptic ruled surface with $e=-1$}
\end{figure}

\noindent {\bf Proof of Theorem \ref{thm:elliptic}(3).} If $b \geq
p+2$, then $L$ satisfies property $N_p$ by Theorem
\ref{thm:elliptic}.(2). For the converse, we prove a stronger statement:\\
\begin{enumerate}
\item[$(*)$] If $1 \leq a \leq p+3$ and $b=p+3-a$, then the line
bundle $N \in \mbox{Pic}X$ in the numerical class $aC_0 + bf$ does
not
satisfy property $N_{p+1}$.\\
\end{enumerate}
From the short exact sequence
\begin{equation*}
0 \rightarrow \mathcal{O}_X (-C_0) \rightarrow \mathcal{O}_X
\rightarrow \mathcal{O}_C \rightarrow 0,
\end{equation*}
we obtain
\begin{eqnarray*}
\cdots  \rightarrow H^1 (X,\wedge^{p+2} \mathcal{M}_N \otimes N)
& \rightarrow & H^1 (C,\wedge^{p+2} \mathcal{M}_N \otimes N_C) \\
& \rightarrow & H^2 (X,\wedge^{p+2} \mathcal{M}_N \otimes
\mathcal{O}_X (-C_0)\otimes N)  \rightarrow \cdots
\end{eqnarray*}
where $N_C$ is the restriction of $N$ to $C_0$. Note that\\
\begin{enumerate}
\item[$(i)$] Since $\mbox{deg}(N_C)=a+b=p+3$, $N_C$ fails to
satisfy property $N_{p+1}$ and hence $H^1 (C,\wedge^{p+2}
\mathcal{M}_N \otimes N_C)\neq 0$ by Lemma \ref{lem:degenrate}.
\item[$(ii)$] The short exact sequence
\begin{equation*}
0 \rightarrow \mathcal{M}_N \rightarrow H^0 (X,N) \otimes
\mathcal{O}_X \rightarrow N \rightarrow 0,
\end{equation*}
gives the following long exact sequence:
\begin{eqnarray*}
\cdots \rightarrow H^1 (X,\wedge^{p+1} \mathcal{M}_N \otimes
\mathcal{O}_X (-C_0)\otimes N^2)
 \rightarrow H^2 (X,\wedge^{p+2} \mathcal{M}_N \otimes
\mathcal{O}_X (-C_0)\otimes N)\\ \rightarrow \wedge^{p+2}H^0
(X,N)\otimes H^2 (X, \mathcal{O}_X (-C_0)\otimes N) \rightarrow
\cdots.
\end{eqnarray*}
It is clear that $H^2 (X, \mathcal{O}_X (-C_0)\otimes N)=0$. Also
by Lemma \ref{lem:generalvanishing},
\begin{equation*}
H^1 (X,\wedge^{p+1} \mathcal{M}_N \otimes \mathcal{O}_X
(-C_0)\otimes N^2)=0
\end{equation*}
and hence $H^2 (X,\wedge^{p+2} \mathcal{M}_N \otimes \mathcal{O}_X
(-C_0)\otimes N)=0$.\\
\end{enumerate}
By $(i)$ and $(ii)$, $H^1 (X,\wedge^{p+2} \mathcal{M}_N \otimes
N)\neq0$ which completes the proof. \qed \\

\noindent Results on higher syzygies of elliptic ruled surfaces
with $e=-1$ are summarized in Figure.1. In this figure, the
integral points of the coordinate plane represent the classes of
$\mbox{Num}(X)$. Theorem \ref{thm:GP} and Theorem
\ref{thm:elliptic} imply that the line bundles contained in the
upper right area of the straight line or the higher crooked line
should satisfy property $N_p$. The above conjecture by F. J.
Gallego and B. P. Purnaprajna suggest that the line bundles
contained in the right upper regions of the lower crooked line
should satisfy property $N_p$.\\

\noindent {\bf Proof of Theorem \ref{thm:main}.} From the short
exact sequence
\begin{equation*}
0 \rightarrow \mathcal{O}_X (-C_0) \rightarrow \mathcal{O}_X
\rightarrow \mathcal{O}_C \rightarrow 0,
\end{equation*}
we obtain
\begin{eqnarray*}
\cdots \rightarrow & H^1 (X,\wedge^j \mathcal{M}_L \otimes
\mathcal{O}_X (-C_0)\otimes L) \rightarrow H^1 (X,\wedge^j
\mathcal{M}_L \otimes L) \\
\rightarrow & H^1 (C,\wedge^j \mathcal{M}_L \otimes L_C)
\rightarrow  H^2 (X,\wedge^j \mathcal{M}_L \otimes \mathcal{O}_X
(-C_0)\otimes L)  \rightarrow \cdots.
\end{eqnarray*}\\
Note that
\begin{enumerate}
\item[$(i)$] Since $\mbox{deg}(L_C)=b-ae \geq 2g+1+p$, $H^1
(C,\wedge^j \mathcal{M}_L \otimes L \otimes \mathcal{O}_C)=0$ for
all $1 \leq j \leq p+1$ by Lemma \ref{lem:degenrate}.
\item[$(ii)$] By Lemma \ref{lem:generalvanishing}, $H^1
(X,\wedge^j \mathcal{M}_L \otimes \mathcal{O}_X (-C_0) \otimes
L)=0$ for $1 \leq j \leq p+1$ if
\begin{equation*}
a-1 \geq p+1~~~~\mbox{and}~~~~b-(a-1)e >
\frac{(p+1)(b-ae)}{b-ae-g}+2g-2.
\end{equation*}
Let $\nu = b-ae$. Then the second inequality is equivalent to
\begin{equation*}
\nu^2 - (3g-e-1+p)\nu+(2g^2-2g-eg)>0
\end{equation*}
and it is a tedious calculation that this inequality holds if $0
\leq e \leq g-3$ and $b-ae\geq 3g-1-e+p$ or if $e \geq g-2$ and
$b-ae\geq 2g+1+p$.
\item[$(iii)$] By the same way as in $(ii)$, when $e \geq g-1$,
\begin{equation*}
H^1 (X,\wedge^{p+2} \mathcal{M}_L \otimes \mathcal{O}_X (-C_0)
\otimes L)=0
\end{equation*}
if $a \geq p+3$ and $b -ae \geq 2g+1+p$. \item[$(iv)$] By Lemma
\ref{lem:generalvanishing}, $H^2 (X,\wedge^j \mathcal{M}_L \otimes
\mathcal{O}_X (-C_0)\otimes L)=0$ if $j \leq a$.
\end{enumerate}
By $(i)$ and $(ii)$, $H^1 (X,\wedge^j \mathcal{M}_L \otimes L)=0$
for all $1 \leq j \leq p+1$ which implies that $(X,L)$ satisfies
property $N_p$ by Lemma \ref{lem:criterion}. Therefore $(1)$ and
$(2)$ are proved. By $(iii)$ and $(iv)$, $H^1 (X,\wedge^{p+2}
\mathcal{M}_L \otimes L) \cong H^1 (C,\wedge^{p+2} \mathcal{M}_L
\otimes L_C)$ which
implies $(3)$. \qed\\

\noindent {\bf Proof of Corollary \ref{cor:amplecase}.} Put $A_i =
a_i C_0 + b_i f$. Since $e \geq 0$, $A_i$ is ample if and only if
$a_ \geq 1$ and $b_i -a_i e \geq 1$. Let $L= K_X + A_1 + \cdots
+A_q$ which is in the numerical class of $(\sum a_i -2 ) C_0 +
(\sum b_i +2g-2-e)f$. Note that
\begin{eqnarray*}
(\sum b_i +2g-2-e)- e(\sum a_i -2 ) & = & \sum (b_i - e a_i )+2g-2 +e \\
                                    & \geq &  q+2g-2+e.
\end{eqnarray*}
To apply Theorem  \ref{thm:main}, it is necessary that $\sum a_i
-2 \geq p+2$ which is true since we assume that $q
\geq 4+p$. Therefore $L$ satisfies property $N_p$ if\\\\
$(1)$ $q+2g-2+e \geq 3g-1-e+p$ when $0 \leq e \leq g-3$, and \\
$(2)$ $q+2g-2+e \geq 2g+1+p$ when $e \geq g-2$\\\\
or equivalently\\\\
$(1)'$ $q \geq \mbox{max} \{ g+1-2e, 4 \} +p$ when $0 \leq e \leq g-3$, and\\
$(2)'$ $q \geq 3-e+p$ when $e \geq g-2$.\\\\
Combining the assumption $q \geq 4+p$ and $(1)'$, $(2)'$, we get
the desired statement.
\end{proof}

\section{Open questions and Conjectures}
\noindent For a given normally generated very ample line bundle,
it is one of the most natural questions that for which $p$ the
line bundle satisfies property $N_p$. Along this line, we suggest
open questions related to higher syzygies of line bundles on ruled
surfaces. Let $X$ be a ruled surface over a smooth projective
curve $C$ of genus $g$. Let $L \in \mbox{Pic}X$ be a line bundle
in the numerical class of $aC_0 +bf$ with $a\geq1$ and let $L_C$
be the restriction of $L$ to $C_0$.

\subsection{Semistable or Unstable Cases} When $e \geq 0$, the followings are
known:\\\\
$(1)$ (Theorem 1.2, \cite{ES}) If $b-ae \geq 3g-1+p$, then
property $N_p$ holds for $L$.

\noindent $(2)$ (Theorem \ref{thm:main}.(1)) When $0 \leq e \leq
g-2$, $L$ satisfies property $N_p$ if $a \geq p+2$ and

$b-ae\geq 3g-1-e+p$.

\noindent $(3)$ (Theorem \ref{thm:main}.(2)) When $e \geq g-1$,
$L$ satisfies property $N_p$ if $a \geq p+2$ and

$b-ae\geq 2g+1+p$.\\

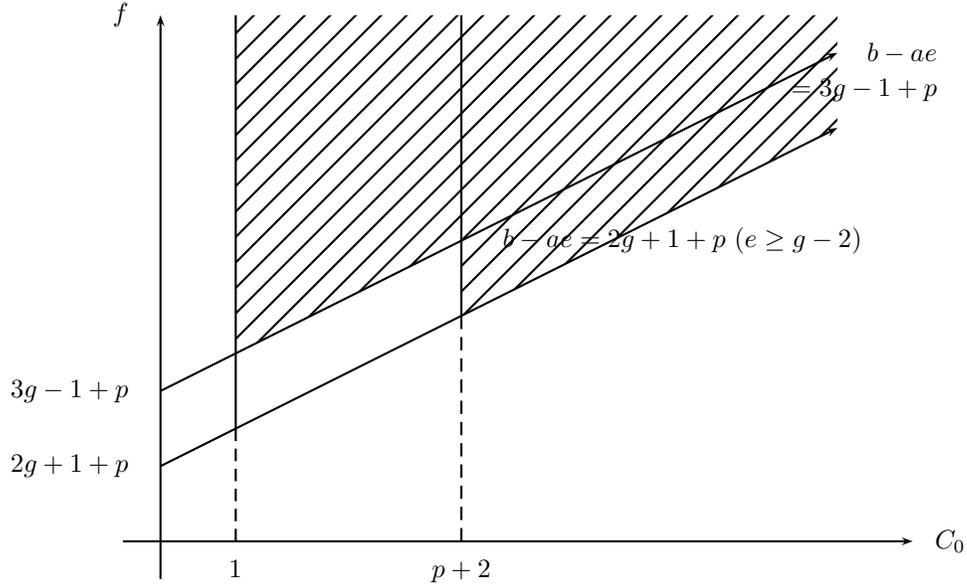
\begin{figure}[t]
\begin{center}
\begin{pspicture}(-1,0)(11,8)
\footnotesize \psset{xunit=.5cm,yunit=.5cm}

\psaxes[labels=none,ticks=none]{->}(0,0)(-1,-1)(20,14)

\uput[180](-.5,2){$2g+1+p$} \uput[180](-.5,4){$3g-1+p$}
\uput[180](-.5,14){$f$} \uput[180](19,8){$b-ae=2g+1+p~(e \geq
g-2)$} \uput[180](21,13){$b-ae$} \uput[180](21,12){$=3g-1+p$}

\rput[t](2,-.5){1} \rput[t](8,-.5){$p+2$} \rput(21,0){$C_0$}

\psline[linestyle=dashed](2,0)(2,3) \psline(2,3)(2,14)
\psline[linestyle=dashed](8,0)(8,6) \psline(8,6)(8,14)
\pspolygon[linestyle=none,fillstyle=hlines,hatchsep=6pt](2,14)(2,5)(8,8)(8,6)(18,11)(18,14)

\psline[arrows=->](0,2)(18,11) \psline[arrows=->](0,4)(18,13)

\end{pspicture}
\end{center}
\caption{Ruled surface with $e \geq 0$}
\end{figure}

\noindent These facts are summarized in Figure.1. In this figure,
the integral points of the coordinate plane represent the classes
of $\mbox{Num}(X)$ and the shadowed regions contain the line
bundles which satisfy property $N_p$. From these results and
Green's ``$2g+1+p$" theorem, the most natural hope in general would be the following:\\

\noindent{\bf Conjecture 1.} $L$ satisfies
property $N_p$ if $b-ae\geq 2g+1+p$.\\

\noindent This conjecture implies Mukai's conjecture for ruled
surfaces with $e \geq 0$ as in the case of Corollary
\ref{cor:amplecase}.(2). More precisely, this guarantees the following: \\

\noindent $(*)$ Let $L=K_X + A_1 + \cdots +A_q$ where $A_i$ is
ample and $q \geq 3$. Then $L$ satisfies

property $N_{p}$ if $q \geq 3-e+p$.\\

\noindent Also in the direction of Theorem \ref{thm:main}.(3), we make the following:\\

\noindent{\bf Conjecture 2.} When $b-ae=2g+p$, $L$ fails to
satisfy property $N_p$ if and only if

either $(i)$  $C$ is hyperelliptic or

$(ii)$ $L_C =K_C + D$ for some effective divisor $D$ of degree
$p+2$.\\

\noindent These two conjectures imply that when $e\geq0$, the
higher syzygies of line bundles would be closely connected with
that of the minimal section. Theorem \ref{thm:elliptic} implies
exactly that these are true for elliptic ruled surfaces.

\subsection{Stable Cases} When $e < 0$, it is known that\\

\noindent $(*)$ (Theorem 1.2, \cite{ES}) If $2b-ae \geq 6g-2+2p$,
then Property $N_p$ holds for $L$.\\

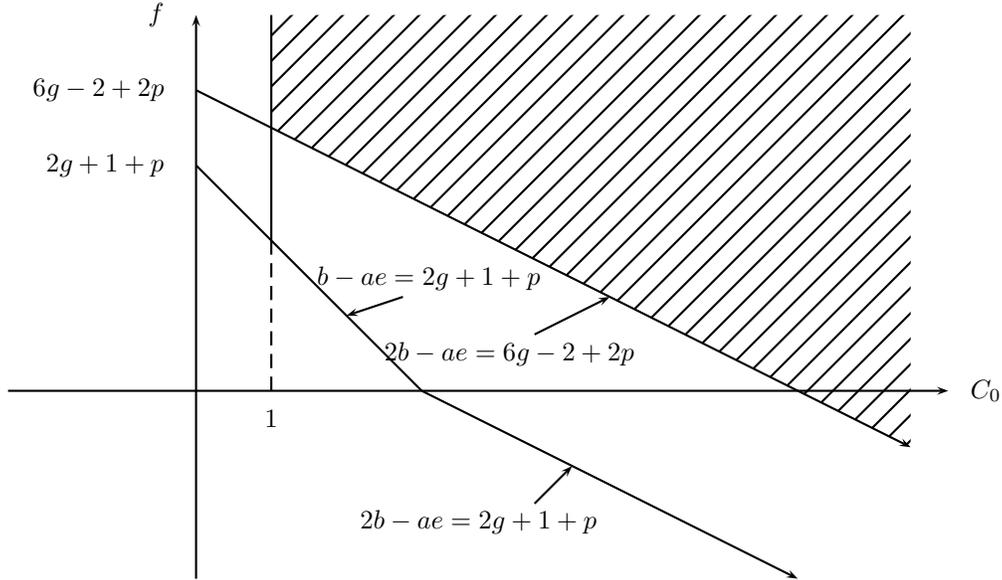
\begin{figure}[t]
\begin{center}
\begin{pspicture}(-2,-2)(11,6)
\footnotesize \psset{xunit=.5cm,yunit=.5cm}

\psaxes[labels=none,ticks=none]{->}(0,0)(-5,-5)(20,10)

\uput[180](-.5,6){$2g+1+p$} \uput[180](-.5,8){$6g-2+2p$}
\uput[180](-.5,10){$f$} \uput[180](12,1){$2b-ae=6g-2+2p$}
\uput[180](9.5,3){$b-ae=2g+1+p$}
\uput[180](11,-3.5){$2b-ae=2g+1+p$}

\rput[t](2,-.5){1} \rput(21,0){$C_0$}

\psline[linestyle=dashed](2,0)(2,4) \psline(2,4)(2,10)
\pspolygon[linestyle=none,fillstyle=hlines,hatchsep=6pt](2,10)(2,7)(19,-1.5)(19,10)

\psline[arrows=->](0,8)(19,-1.5) \psline[arrows=->](5.5,2.5)(4,2)
\psline[arrows=->](9,1.5)(11,2.5) \psline[arrows=->](9,-3)(10,-2)
\psline(0,6)(6,0) \psline[arrows=->](6,0)(16,-5)

\end{pspicture}
\end{center}
\caption{Ruled surface with $e < 0$}
\end{figure}

\noindent Note that when $e <0$, there may be other extremal
curves as discussed in $\S 3$. Indeed results of Yuko
Homma\cite{Homma2} and Gallego-Purnaprajna\cite{GP2}
about elliptic ruled surfaces with $e=-1$ make the following conjecture affirmative.\\

\noindent {\bf Conjecture 3.}  $L$ satisfies Property $N_p$
if $b-ae \geq 2g+1+p$ and $2b-ae \geq 2g+1+p$. \\

\noindent We remark that this is proved for $g=1$ and $a=1$, i.e.,
elliptic surface scrolls. As in the case $e \geq 0$, this
guarantees Mukai's conjecture for ruled surfaces $e<0$. In Figure
2, the integral points of the coordinate plane represent the
classes of $\mbox{Num}(X)$ and the shadowed regions contain
the line bundles which satisfy property $N_p$.\\

\end{document}